\documentclass{amsart}
\usepackage{amssymb,latexsym}
\usepackage{graphicx}
\newtheorem{theorem}{Theorem}
\newtheorem{lemma}{Lemma}
\newtheorem{definition}{Definition}
\newtheorem{notation}{Notation}
\newtheorem{example}{Example}
\newtheorem{corollary}{Corollary}

\begin{document}
\title{ FLUCTUATIONS OF MULTI-DIMENSIONAL\\ KINGMAN-L\'EVY PROCESSES }

\author{Thu Nguyen}
\address{Department of Mathematics;
         International University, HCM City;
         No.6 Linh Trung ward, Thu Duc District, HCM City;
         Email: nvthu@hcmiu.edu.vn}
\date{August 10, 2009}
\begin{abstract} In the recent paper \cite{Ng5} we have introduced a method of studying the multi-dimensional Kingman convolutions and their associated stochastic processes by embedding them
into some multi-dimensional ordinary convolutions which allows to study multi-dimensional Bessel processes in terms of the cooresponding Brownian motions. Our further aim in this paper is to
introduce  k-dimensional Kingman-L\'evy (KL) processes and prove some of their fluctuation
properties which  are analoguous to that of k-symmetric L\'evy processes. In particular, the L\'evy-It\^o decomposition and  the  series representation of Rosi\'nski type for k-dimensional KL-processes are obtained. 
    \end{abstract} 
        \maketitle{Keywords and phrases: Cartesian products of Kingman convolutions; Rayleigh distributions}
\section{Introduction. notations and prelimilaries}
 The purpose of this paper is to introduce and study the multivariate KL processes defined in terms of multicariate Kingman convolutions. To begin with we review the following information of the Kingman convolutions and their Cartesian products.

  Let $\mathcal P:=\mathcal P(\mathbb R^+)$  denote the set of all probability measures (p.m.'s) on the positive half-line   $\mathbb R^+$. Put, for each continuous bounded function f on $\mathbb R^{+}$,
 \begin{multline}\label{astKi}
\int_{0}^{\infty}f(x)\mu\ast_{1,\delta}\nu(dx)=\frac{\Gamma(s+1)}{\sqrt{\pi}\Gamma(s+\frac{1}{2})}\\
\int_{0}^{\infty}\int_{0}^{\infty}\int_{-1}^{1}f((x^2+2uxy+y^2)^{1/2})(1-u^2)^{s-1/2}\mu(dx)\nu(dy)du,
\end{multline}
  where $\mu\mbox{ and }\nu\in\mathcal P\mbox{ and }\delta=2(s+1)\geq1$ (cf. Kingman\cite{Ki} and
Urbanik\cite{U1}). The convolution algebra $(\mathcal{P},\ast_{1,\delta})$ is
 the most important example of Urbanik convolution algebras (cf. Urbanik\cite{U1}). In language of the
  Urbanik convolution algebras, the {\it characteristic measure}, say $\sigma_s$, of the Kingman convolution
  has the Rayleigh density
  \begin{equation}\label{Ray}
  d\sigma_s(y)= \frac{2{(s+1)^{s+1}}}{\Gamma(s+1)}y^{2s+1}\exp{(-(s+1)y^2)}dy
  \end{equation}
with the characteristic exponent $\varkappa=2$ and the kernel
$\Lambda_s$
\begin{equation}\label{eq:Lam}
 \Lambda_s(x)= \Gamma(s+1) J_{s}(x)/(1/2x)^{s},
\end{equation}
where $J_s(x)$ denotes the Bessel function of the first kind,
\begin{equation}\label{eq:Bessel}
J_s(x):= \Sigma_{k=0}^{\infty} \frac{(-1)^k
(x/2)^{\nu+2k}}{k!\Gamma(\nu+k+1)}.
\end{equation}
  It is known (cf. Kingman \cite{Ki}, Theorem 1), that the kernel $\Lambda_s$ itself is an
ordinary characteristic function (ch.f.) of a symmetric p.m., say $F_s$, defined on the
interval [-1,1]. Thus, if $\theta_s$ denotes a random variable (r.v.) with distribution
$F_s$ then for each $t\in \mathbb R^+$,
\begin{equation}\label{eq:LamThe}
\Lambda_s(t)= E\exp{(it\theta_s)}=\int_{-1}^1\cos{(tx)}dF_s(x).\end{equation}
 Suppose that $X$ is a nonnegative r.v. with distribution $\mu\in\mathcal{P}$
 and $X$ is independent of $\theta_s$.  The {\it radial characteristic function}
 (rad.ch.f.) of $\mu$, denoted by $\hat\mu(t),$ is defined by
  \begin{equation}\label{ra.ch.f.}
\hat\mu(t) = E\exp{(itX\theta_s)} = \int_0^{\infty}
\Lambda_s(tx)\mu(dx),
\end{equation}
 for every $t\in \mathbb R^{+}$.
 The characteristic measure of the Kingman convolution $\ast_{1, \delta}$, denoted by $\sigma_s$,
has the Maxwell density function
\begin{equation}\label{Maxwell density}
\frac{d\sigma_s(x)}{dx}=\frac{2(s+1)^{s+1}}{\Gamma(s+1)}x^{2s+1}exp\{-(s+1)x^2\}, \quad(0<x<\infty).
\end{equation}
and the rad.ch.f.
\begin{equation}
\hat\sigma_s(t)=exp\{-t^2/4(s+1)\}.
\end{equation}

Let $\tilde P:=\tilde{\mathcal P}(\mathbb R)$ denote the class of symmetric p.m.'s on $\mathbb R.$ Putting, for every $G\in \mathcal P$,
\begin{equation*}
 F_s(G)=\int_0^{\infty}F_{cs} G(dc),
 \end{equation*}
  we get a continuous homeomorphism from the Kingman convolution algebra $(\mathcal{P},\ast_{1,\delta})$ onto the ordinary convolution algebra $(\tilde{\mathcal P}, \ast)$ such that
\begin{eqnarray}\label{homeomorphism1}
 F_s\{G_1\ast_{1, \delta}G_2\}&=&(F_sG_1)\ast(F_sG_2) \qquad (G_1, G_2\in \mathcal P)\\
 F_s\sigma_s&=&N(0, 2s+1)
 \end{eqnarray}
 which shows that every Kingman convolution can be embedded into the ordinary convolution $\ast$.
 
   Denote by $ \mathbb {R}^{+k}, k=1,2,...$ the k-dimensional nonnegative cone
 of $ \mathbb {R}^{k}$ and $\mathcal{P}(\mathbb {\mathbb R}^{+k})$ the class of
all p.m.'s on $\mathbb  {\mathbb R}^{+k}$ equipped with the weak convergence. In
the sequel, we will denote the multidimensional vectors and random vectors (r.vec.'s)
and their distributions by bold face letters.

For each point z of any set $A$ let $\delta_z$ denote the Dirac measure (the unit mass) at
the point z. In particular, if $\mathbf  x=(x_1, x_2,\cdots,x_k)\in
 \mathbb R^{k+}$, then
\begin{equation}\label{proddelta}
\delta_{\mathbf {x}} =
\delta_{x_1}\times\delta_{x_2}\times \ldots\times\delta_{x_k},\quad (k\; times),
\end{equation}
where the sign $"\times"$ denotes the Cartesian product of
measures.
  We put, for $\mathbf {x} = (x_1,\cdots, x_k)\mbox{ and }\mathbf {y} =
(y_1,y_2,\cdots, y_k)\in \mathbb R^{+k},$
  \begin{equation}\label{convdeltas}\delta_{\mathbf x}\bigcirc_{s, k} \delta_{\mathbf {y}} = \{\delta_{x_1}\circ _s \delta_{y_1}\} \times\{\delta_{x_2} \circ _s\delta_{y_2}\}
\times\cdots\
 \times \{\delta_{x_k} \circ_s \delta_{y_k}\},\quad (k\; times),
\end{equation}
here and somewhere below for the sake of simplicity we denote the
Kingman convolution operation $\ast_{1,\delta}, \delta=2(s+1)\ge 1$ simply by $\circ_{s}, s\ge \frac{!}{2}.$
    Since convex combinations of p.m.'s of the form
(\ref{proddelta}) are dense in $\mathcal P(\mathbb R^{+k})$ the
relation (\ref{convdeltas}) can be extended to arbitrary p.m.'s $
\mathbf{G}_1 \mbox{  and } \mathbf{G}_2\in\mathcal{P}( \mathbb R^{+k})$.
Namely, we put
\begin{equation}\label{convF}
\mathbf {G}_1 \bigcirc_{s, k} \mathbf {G}_2 = \iint\limits_{ \mathbb R^{+k}}
  \delta_{\mathbf {x}} \bigcirc_{s, k} \delta_{\mathbf {y}} {\mathbf G}_1(d\mathbf {x}) {\mathbf G}_2(d\mathbf {y})
 \end{equation} which means that for each continuous bounded function $\phi$ defined on $\mathbb R^{+k}$
 \begin{equation}\label{convof}
\int\limits_{\mathbb R^{+k}} \phi({\mathbf z}) {\mathbf G}_1 \bigcirc_{s, k} {\mathbf G}_2  (d{\mathbf z})= \iint\limits_{ \mathbb R^{+k}}\big\{\int\limits_{\mathbb R^{+k}} \phi({\mathbf z})  \delta_{{\mathbf x}} \bigcirc_{s, k} \delta_{{\mathbf y}}(d{\mathbf z})\big\}{ \mathbf G}_1(d{\mathbf x}) {\mathbf G}_2(d{\mathbf y}).
\end{equation}
 In the sequel, the binary
operation $\bigcirc_{s, k}$ will be called {\it the k-times Cartesian
product of Kingman convolutions} and the pair $(\mathcal P( \mathbb R^{+k}), \bigcirc_{s, k})$ will be called
{\it the k-dimensional Kingman convolution algebra}. It is easy to show that the
binary operation $\bigcirc_{s, k}$ is continuous in the weak topology
which together with (\ref{astKi}) and (\ref{convF}) implies the
following theorem.
 \begin{theorem}\label{Theo:Kingmanalgebra} The pair $(\mathcal P{( \mathbb R^{+k})} ,\bigcirc_{s,  k})$
 is a commutative
topological semigroup with $\delta_{\mathbf  0}$ as the unit element.
Moreover, the operation $\bigcirc_{s, k}$ is distributive w.r.t.
convex combinations of p.m.'s in $\mathcal P( \mathbb R^{+k})$.
\end{theorem}
 \  For every ${\mathbf  G}\in\mathcal P( \mathbb R^{+k})$ the
k-dimensional rad.ch.f. $\hat{{\mathbf  G}}({\mathbf  t}), {\mathbf
t}=(t_1, t_2, \cdots t_k)\in \mathbb R^{+k},$ is defined by
\begin{equation}\label{k-ra.ch.f.}
\hat{\mathbf  G}(\mathbf  {t})=\int\limits_{\mathbb R^{+k}}
 \prod_{j=1}^{k}\Lambda_s(t_jx_j){ \mathbf  G}(\mathbf {dx}),
 \end{equation}
 where $\mathbf {x}=(x_1, x_2, \cdots x_k)\in  \mathbb R^{+k}.$
 Let $\mathbf{\Theta_s} = \{\theta_{s, 1},\theta_{s, 2}, \cdots ,\theta_{s, k}\}$, where $\theta_{s, j}$ are independent r.v.'s with the same distribution
 $F_s $.
 Next, assume that $ {\mathbf X}=\{X_1, X_2,..., X_k\}$ is a k-dimensional
r.vec. with distribution $\mathbf{G}$ and $\mathbf{X}$ is independent of
$\mathbf{\Theta}_s$. We put
     \begin{equation}\label{[Theta,X]}
[{\mathbf\Theta}_s,{\mathbf X}]=\{{\theta_{s,1} X_1, \theta_{s, 2} X_2,...,\theta_{s, k}X_k}\}.
\end{equation}
 Then, the following formula is equivalent to (\ref{k-ra.ch.f.}) (cf. \cite{Ng4})
 \begin{equation}\label{multiradchf}
 \widehat{\mathbf G}({\mathbf t})=Ee^{i<{\mathbf t},[{\mathbf\Theta_s, \mathbf X}]>},\qquad ({\mathbf t}\in \mathbb R^{+k}).
 \end{equation}
The Reader is referred to Corollary 2.1, Theorems 2.3 \& 2.4 \cite{Ng4} for the  principal properties of  the above rad.ch.f.
 Given $s\ge -1/2$  define a map $F_{s, k}:  \mathcal  P(\mathbb R^{+k}) \rightarrow  \mathcal P(\mathbb R^k)$ by
 \begin{equation}\label{k-map}
 F_{s, k}({\mathbf G})=\int\limits_{\mathbb R^{+k}} (T_{c_1}F_s)\times(T_{c_2}F_s)\times  \ldots\times(T_{c_k}F_s) {\mathbf G}(d{\mathbf  c}),
 \end{equation}
 here and in the sequel, for a distribution  $\mathbf G$ of a r.vec. $\mathbf X$ and a real number c we denote by $T_c{\mathbf G}$ the distribution of $c{\mathbf X}$.
  Let us denote by $\tilde{ \mathcal P}_{s, k}(\mathbb{R}^{+k})$ the sub-class of  $\mathcal P(\mathbb R^k)$ consisted of all p.m.'s defined by the right-hand side of (\ref{k-map}).
 By virtue of (\ref{k-ra.ch.f.})-(\ref{k-map}) one can prove the following theorem.
 \begin{theorem}\label{symmconvo}
The set $\tilde{ \mathcal P}_{s, k}(\mathbb{R}^{+k})$ is closed w.r.t. the weak convergence and the ordinary convolution $\big.\ast$ and the following equation holds
\begin{equation}\label{Fourier=rad.ch.f.}
\hat{\mathbf G}({\mathbf t})=\mathcal F(F_{s, k}({\mathbf G}))({\mathbf t})\qquad ({\mathbf t}\in {\mathbb R^{+k}})
\end{equation}
where $\mathcal F({\mathbf K})$ denotes the ordinary characteristic function (Fourier transform)
of a p.m. ${\mathbf K}$. Therefore, for any ${\mathbf G}_1\mbox{ and } {\mathbf G}_2\in \mathbb R^{+k}$
\begin{equation}\label{convolequality}
F_{s, k}({\mathbf G}_1)\big.\ast F_{s, k}({\mathbf G}_2)=F_{s, k}({\mathbf G}_1\bigcirc_{s, k}{\mathbf G}_2)
\end{equation}
and the map $F_{s, k}$ commutes with convex combinations of distributions and scale changes
$T_c, c>0$. Moreover,
\begin{equation}\label{Gaussian-Rayleigh}
F_{s, k}({\Sigma_{s, k}})=N({\mathbf 0}, 2(s+1){\mathbf I})
\end{equation}
where $\Sigma_{s, k}$ denotes the k-dimensional Rayleigh distribution and $N({\mathbf 0}, 2(s+1){\mathbf I})
$ is the symmetric normal distribution on $\mathbb R^k \mbox{ with variance operator } {\mathbf R}= 2(s+1)
{\mathbf I}, {\mathbf I}$ being the identity operator. 
Consequently, Every Kingman convolution algebra $\big( \mathcal P(\mathbb R^{+k}), \bigcirc_{s, k}\big)$ is
embedded in the ordinary convolution algebra $\big( \mathcal P_{s, k}(\mathbb{R}^{+k}), \big.\star \big)$ and the map $F_{s, k}$ stands for a homeomorphism.
\end{theorem}
 
 Let us denote by $\mathcal E=\{{\mathbf e}=(e_1, e_2, \ldots, e_k), e_j=\pm , j=1, 2, \ldots, k \}$.
  It is convenient to regard the elements of $\mathcal E$ as sign vectors. Denote  $\mathbb R^{+k}_{\mathbf e} =\{[{\mathbf e},{\mathbf x} ]: {\mathbf x}\in \mathbb R^{+k}\}, \mbox{ where } [{\mathbf e},{\mathbf x} ]:=(e_1x_1, e_2x_2, \ldots, e_kx_k).$ Then  the family             $\{\mathbb R^{+k}_{\mathbf e}, {\mathbf e} \in\mathcal E \}$  is a partition of the space $\mathbb R^k.$ If $\mathbf X$ is a k-dimensional r.vec. with distribution $\mathbf G,$ the k-symmetrization of $\mathbf G$, denoted by $\tilde{\mathbf G},$ is defined by
  \begin{equation}
  \tilde{\mathbf G}=\frac{1}{2^k} \sum_{\mathbf e \in \mathcal E} S_{\mathbf e} {\mathbf G},
  \end{equation}
  where the operator $S_{\mathbf e}$ is defined by
  \begin{equation}
  S_{\mathbf e}({\mathbf x})=[{\mathbf e},{\mathbf x}] \qquad {\mathbf x}\in{ \mathbb  R^k}
  \end{equation}
  and the symbol $S_{\mathbf e} \tilde{\mathbf G}$ denotes the image of $\mathbf G$ under $S_{\mathbf e}$.
\begin{definition} \label{k-symm}
 We say that a distribution $\mathbf G\in \mathcal P(\mathbb R^k)$ is k-symmetric, if the equation
 $\mathbf G=\tilde{\mathbf G}$ holds.
 \end{definition}

   \begin{definition}\label{k-ID}
 A p.m. ${\mathbf F} \in  \mathcal P(\mathbb R^{+k})$ is called $\bigcirc_{s, k}-$infinitely divisible
 ($\bigcirc_{s, k}-$ID), if for every m=1, 2, \ldots there exists $\mathbf F_m\in  \mathbf P(\mathbb R^{+k})$ such that
 \begin{equation}\label{kID}
{ \mathbf F}={\mathbf F}_m\bigcirc_{s, k} {\mathbf F}_m\bigcirc_{s, k}\ldots \bigcirc_{s, k}{\mathbf F}_m\quad (m\;times).
 \end{equation}
 \end{definition}
 \begin{definition}\label{stability}
 $\mathbf F$ is called stable, if for any positive
 numbers a and b there exists a positive number c such that
 \begin{equation}\label{k-stability}
  T_a{\mathbf F}\;{\bigcirc_{s, k}}\;T_b{\mathbf F}=T_c{\mathbf F}
 \end{equation}
\end{definition}
 By virtue of Theorem \ref{symmconvo} it follows that the following theorem holds.
 \begin{theorem}\label{equivdef}
 A p.m. $\mathbf G\mbox{ is } \bigcirc_{s, k}-ID$, resp. stable if and only if
 $H_{s, k}({\mathbf G})$ is ID, resp. stable,  in the usual sense.
 \end{theorem}
   The following theorem gives a representation of rad.ch.f.'s of  $\bigcirc_{s, k}-$ID distributions. The proof is  a verbatim reprint of that for  (\cite{Ng4}, Theorem 2.6).
 \begin{theorem}\label{LevyID} A p.m. $\mu\in ID(\bigcirc_{s, k})$ if and only if
there exist a $\sigma$-finite measure M (a L\'evy's measure) on
  $ \mathbb R^{+k}$ with the property that $M({\mathbf 0})=0,  {\mathbf M}$  is finite outside every neighborhood of ${\mathbf 0}$ and
\begin{equation}\label{integrable w. r. t. weight function}
\int_{\mathbb R^{+k}}\frac {\|{\mathbf x}\|^2} {1+\|{\mathbf x}\|^2}
{\mathbf M}(d{\mathbf x}) < \infty
\end{equation}
 and for each ${\mathbf t}=(t_1,...,t_k)\in
\mathbb R^{+k}$
\begin{equation}\label{Levy-Kintchine for k-dim.rad. ch. f.}
 -\log{\hat{\mu}({\mathbf  t})}=\int_{\mathbb R^{+k}}\{1-\prod_{j=1}^{k}\Lambda_s(t_jx_j)\} \frac
{1+\|{\mathbf x}\|^2} {\|{\mathbf x}\|^2} M({\mathbf {dx}}),
\end{equation}
where, at the origin $\mathbf{0}$, the integrand on the right-hand
side of (\ref{Levy-Kintchine for k-dim.rad. ch. f.}) is assumed to
be
\begin{equation}\label{limiting integrand}
lim_{\|\mathbf
{x}\|\rightarrow 0 }\{1-\prod_{j=1}^k
 \Lambda_s(t_jx_j)\} \frac {1+\|\mathbf x\|^2}
{\|\mathbf {x}\|^2}=
\Sigma_{j=1}^k \lambda^2_j t_j^2   
\end{equation}
for nonnegative $\lambda_j, j=1, 2,...,k.$
 In particular, if $ M=0, \mbox{ then } \mu $ becomes a Rayleighian distribution with
the rad.ch.f (see definition \ref{Rayleigh})
\begin{equation}\label{kRayleighian rad. ch. f.}
-\log{\hat{\mu}({\mathbf t})}=\frac{1}{2}\sum_{j=1}^k \lambda^2_j
t_j^2,\quad {\mathbf t}\in \mathbb R^{+k},
\end{equation}
 for some nonnegative $\lambda_j, j=1,...,k.$
 Moreover, the representation (\ref{Levy-Kintchine for k-dim.rad. ch.
 f.}) is unique.
  \end{theorem}
   An immediate consequence of the above theorem is the following:
\begin{corollary}\label{Cor:Pair}
Each distribution $\mu\in ID(\bigcirc_{s, k})$ is uniquely determined by the pair $[\mathbf{M}, \pmb {\lambda}]$, where $\mathbf{M}$ is a  Levy's measure
on $\mathbb R^{+k}$ such that $\mathbf{M}(\mathbf{0})=0,$ $\mathbf{M}$ is finite outsite every neighbourhood of $\mathbf{0}$ and the condition (\ref{integrable w. r. t. weight function})
 is satisfied and $\pmb{\lambda}:=\{\lambda_1, \lambda_2,\cdots \lambda_k\}\in \mathbb R^{+k}$ is a vector of nonnegative numbers appearing in (\ref{kRayleighian rad. ch. f.}).
Consequently, one can write $\mu\equiv[\mathbf{M}, \pmb {\lambda}].$\\ \indent
In particular, if $\mathbf{M}$ is zero measure then $\mu=[\pmb{\lambda}]$ becomes a Rayleighian p.m. on $\mathbb R^{+k}$ as defined as follows:
\end{corollary}
  \begin{definition}\label{Rayleigh}
  A  k-dimensional distribution, say $\pmb{\mathbf \Sigma}_{s, k}$, is called  a {\it Rayleigh distribution}, if   \begin{equation}\label{k-dimension Rayleigh}
\pmb{\mathbf  \Sigma}_{s, k}=\sigma_s\times\sigma_s\times\cdots\times\sigma_s \quad
 (k\;times).
 \end{equation}
 Further, a distribution ${\mathbf F}\in \mathcal P(\mathbb R^{+k})$ is called a {\it Rayleighian distribution} if there exist nonnegative numbers $\lambda_r,
   r=1,2 \cdots k $ such that
\begin{equation}\label{k-dimensional rayleighian}
{ \mathbf F}=\{T_{\lambda_1}\sigma_s\}\times \{T_{\lambda_2}\sigma_s\}
 \times\ldots \times\{T_{\lambda_k}\sigma_s\}.
 \end{equation}
 \end{definition}
 \indent
 It is evident that every Rayleigh distribution is  a Rayleighian distribution. Moreover, every Rayleighian distribution is $\bigcirc_{s, k}-$ID. By virtue of (\ref{Maxwell density} )  and (\ref{k-dimension Rayleigh}) it follows that the k-dimensional Rayleigh density is given by
 \begin{equation}\label{density k-dimension Rayleigh}
 g({\mathbf x})=\Pi_{j=1}^k\frac{2^k(s+1)^{k(s+1)}}{\Gamma^k(s+1)}x_j^{2s+1}exp\{-(s+1)||{\mathbf x}||^2\},
 \end{equation}
 where ${\mathbf x}=(x_1, x_2,\ldots, x_k)\in \mathbb R^{+k}$ and the corresponding rad.ch.f. is given by
 \begin{equation}
 \hat\Sigma_{s, k}({\mathbf t})=Exp(-|{\mathbf t}|^2/4(s+1)),\quad {\mathbf t}\in \mathbb R^{+k}.
 \end{equation}
 Finally, the rad.ch.f. of a Rayleighian distribution $\mathbf  F\mbox{ on } \mathbb R^{+k}$ is given by
 \begin{equation}\label{rad.ch.rayleighian}
 \hat{\mathbf F}({\mathbf t})=Exp(-\frac{1}{2}\sum_{j=1}^k\lambda_j^2t_j^2)
 \end{equation}
 where $\lambda_j, j=1, 2, \ldots, k$ are some nonnegative numbers.
\section{Multivariate Bessel processes}
 \section{Multivrariate Kingman-L\'evy processes and their L\'evy-It\^o decomposition}
 Suppose that $\mu_t, t\ge 0$ is continuous semigroup in $ID(\bigcirc_{s, k})$, that is for any $t,s\ge 0$ 
 \begin{equation}\label{Ksemigroup}
 \mu_t\bigcirc_{s, k}\mu_s=\mu_{t+s}
 \end{equation}
 and $\{\mu_t\}$ is continuous at 0 i.e.  
 \begin{equation*}
 lim_{t\rightarrow 0}\mu_t=\delta_{\mbox{0}}.
 \end{equation*}
 By virtue of Theorem \ref{symmconvo} it follows that $\{\mathcal F_{s, k}(\mu_t)\}$ is an ordinary  continuous convolution semigroup on $\mathbb R^k.$
   Putting, for each $\mathbf{x}\in \mathbb {R}^{k+}$ and for every Borel subset $\mathcal{E}\mbox{ of } \mathbb {R}^{k+}, $
 \begin{equation}\label{transitionProbLevy}
 \mathbf{P}(t, \mathcal{E}, \mathbf{x})=\mu_t\bigcirc_{s, k}\delta_{\mathbf{x}}(\mathcal{E})
 \end{equation}
 and using the rad.ch.f. it follows that the family $\{ \mathbf{P}(t, \mathcal{E}, \mathbf{x}), t\ge 0\}$ satisfies the Chapman-Kolmogorov equation and, consequently,  the formula (\ref{transitionProbLevy}) defines transition probabilities of a $\mathbb R^{k+}-$valued  homogeneous strong Markov Feller process $\{\mathbf{X}^{\mathbf{x}}_t, t\ge 0\}$, say, such that
 it is stochastically continuous and has a cadlag version (compare \cite{Ng1}, Theorem 2.6). 
 \begin{definition}\label{KL-process}
 A $\mathbb R^{k+}$-valued stochastic process $\{\mathbf{X}_t, t\ge 0\}$ is called
 a Kingman-L\'evy process, if  $\mathbf{X}_t=$
 
 (i) $ \mathbf{X}_0=\mathbf{0}\qquad (P.1);$
 
 (ii) There exists  a $\mathbb R^{k+}-$valued  homogeneous strong Markov Feller process  having a cadlag version $\{\mathbf{X}^{\mathbf{x}}_t, t\ge 0\}$ with transition probabilities defined by (\ref{transitionProbLevy}) and  $\mathbf{X}_t=\mathbf{X}^{\mathbf{0}}_t, t\ge 0;$
\end{definition}
\section{Fluctuations of Multidimensional Bessel Processes}
\begin{definition}
Let $(W_t, t\ge 0)$ be a d-dimensional Brownian motion (d=1, 2, \ldots). The Euclidean norm of $(W_t)$, denoted by $B_t, t\ge 0$ is called a Bessel process.  
\end{definition}
It has been proved that Bessel processes inherit the well-known characteristics of Brownian motions: They are independent stationary "increments" processes with continuous sample paths. The term 'increment' is defined as follows:
\begin{definition}
For any  $s>u$ the random variable $| W_s-W_u |$ is called an increments of the Bessel process.
\end{definition}
The following theorem gives a L\'evy-Khinczyn representation of the Bessel process in the sense of the Kingman convolution.
\begin{theorem}
The radial characteristic function $\phi(x)$ of the Bessel process $(B_t)$ is of the form
\begin{equation}
\phi(x)=exp\{-\frac{tx^2}{4(s+1)}\}\qquad x, t \ge0
\end{equation}
where d=2(s+1).
\end{theorem}
Since for any $s>u$ the 'increment' of the Bessel process $(B_t)$ is infinitely divisible in the ordinary convolution $\ast$ we have the following representation of the Fourier transform of $B_{s-u}.$
\begin{equation}
\mathcal F_{B_{s-u}}(x)=exp(-(s-u) \psi(x))
\end{equation} 
where $\psi(x)$ is a symmetric characteristic exponent
\begin{equation}
\psi(x)=\frac{1}{2}\sigma^2+\int_0^{\infty}(1-cos\,xv)\Pi(dv)
\end{equation}
where the measure $\Pi$ satisfies the condition \\begin equation
\begin{equation}
\int_0^{\infty}(min(1,x^2)\Pi(dx)<\infty.
\end{equation}
which implies the following L\'evy-It\^o decomposition.
\begin{theorem}(L\'evy-It\^o decomposition)
There exists a  Brownian motion $X^{(1)}_t$ and a compound Poison process $ X^{(2)}_t $  independent of  $X^{(1)}_t$  such that
\begin{equation}
B_t=||W_t ||\overset{d}{=}X^{(1)}_t + X^{(2)}_t\qquad (t\ge 0).
\end{equation}
\end{theorem}
Before stating the Wienner-Hopf factorization (WHf) theorem for Bessel processes we introduce some concepts and notations. The importance of WHf is that it gives us information of the ascending and descending ladder processes. We begin by recalling that for $\alpha, \beta\ge 0$
the Laplace exponents $\kappa(\alpha, \beta)\mbox { and } \hat \kappa(\alpha, \beta)$ of the ascending ladder process $(\hat L^{-1}, \hat H)$ and the descending ladder process $(\hat L^{-1}, \hat H).$ Further, we define
$$\overset {-}{G}_t=sup\{s<t: \overset{-}{X}_s=X_s\} \mbox { and } \underset{-}{G_t}= sup\{s<t: \underset{-}{X_t}=X_s.$$
 \begin{theorem}\label{Wienner-Hopf}(Wienner-Hopf Factorization)
Let $(B_t, t\ge 0)$ be a Bessel process. Denote by ${\mathbf e}_p$ an independent and exponentially distributed random variable.

The pairs $(\overset{-}{G}_{{\mathbf e}_p}, \overset{-}{X}_{{\mathbf e}_p})\mbox{ and }
({\mathbf e}_p-\overset{-}{G}_{{\mathbf e}_p}, \overset{-}{X}_{{\mathbf e}_p}-X_{{\mathbf e}_p})$
are independent and infinitely divisible, yielding the factorization
\begin{equation}\label{Wienner-Hopf}
\frac{p}{p-i\nu+\psi(\theta)}=\Psi^{+}(\nu, \theta).\Psi^{-}(\nu, \theta)\qquad \nu, \theta\in \mathbb R,
\end{equation}
$\psi^{+}, \psi^{-}$ being Fourier transforms and called the Wienner-Hopf factors.
\end{theorem}

\section{Levy-Ito decomposition of Kingman-Levy processes}
\

\end{document}